\documentclass{amsart}
\newtheorem{theorem}{Theorem}[section]
\newtheorem{prop}[theorem]{Proposition}
\newtheorem{lemma}[theorem]{Lemma}
\newtheorem{corollary}[theorem]{Corollary}
\newtheorem{remark}[theorem]{Remark}
\newtheorem{definition}[theorem]{Definition}
\newcommand{\spin}{$Spin^c$-structure }

\usepackage{amsfonts}
\usepackage{amssymb}
\usepackage{amsmath}
\usepackage{amscd}

\begin{document}
\title[Symplectic Surfaces and Generic $J$-holomorphic Structures]{Symplectic Surfaces and Generic $J$-holomorphic Structures on 4-Manifolds}
\author{Stanislav Jabuka}
\address{Department of Mathematics, Columbia University, New York, NY 10027, USA}
\email{jabuka@math.columbia.edu}
\thanks{The author was partially supported by the VIGRE postdoctoral program of the NSF} 
\begin{abstract}
It is a well known fact that every embedded symplectic surface $\Sigma$ in a symplectic four-manifold 
$(X^4,\omega )$ can be made $J$-holomorphic for some almost-complex structure $J$ compatible with $\omega$. 
In this paper we investigate when such a $J$ can be chosen generically in the sense of Taubes (for definition, see below). The main result is stated in 
Theorem \ref{main} below. As an application we give examples of smooth and non-empty Seiberg-Witten and Gromov-Witten moduli spaces
whose associated invariants are zero. 
\end{abstract}
\subjclass{53D35, 57R17; 53D45, 53D05}
\maketitle
\section{Introduction}

To set up the background for the main theorem below, let $C \subset X$ be a connected, symplectic surface  embedded in the minimal symplectic 4-manifold 
$X$ with symplectic form $\omega$. 
It is a well known fact that $C$ can be 
made $J$-holomorphic for some almost-complex structure $J$ compatible with $\omega$. This paper investigates when $J$ can be chosen from a generic set of almost-complex structures.  We start 
by recalling what {\sl generic} means in our setting.

For a given $E\in H_2(X;\mathbb{Z})$, set 
\begin{equation} \label{indexofd}
d=\frac{1}{2}(E^2-K\cdot E)
\end{equation}
where $K$ is the canonical class of $X$
associated to $\omega$. Introduce  $\mathcal{A} _d (X) $ as the set of pairs $(J,\Omega)$ with $J$ an almost-complex structure compatible with $\omega$ and 
$\Omega$ a set 
of $d$ distinct points of $X$. $\mathcal{A} _d (X) $ has the structure of a 
smooth manifold inherited from the Frechet manifold $C^{\infty}(\mbox{End}(TX)\times \mbox{Sym}^d (X))$. 

Each $J$-holomorphic curve $C$ comes equipped with a linear operator 
$$D_C:C^\infty (N_C) \rightarrow C^\infty (N_C \otimes T^{0,1}C)$$
obtained from the linearization of the 
generalized Cauchy-Riemann operator $\overline{\partial _C}$. Here $N_C$ is the normal bundle of $C$ in $X$.  The operator 
$D_C$ is elliptic and its (complex) index is given by $d$ as defined in \eqref{indexofd} with $E=[C]$.
In the case when $C$ contains all points of $\Omega$, let $ev_\Omega:C^\infty (N_C) \rightarrow \oplus _{p\in \Omega} N_p$ be the evaluation map 
associated to $\Omega$. If $d=0$, we say that $D_C$ is non-degenerate if $\mbox{Coker}(D_C)=\{ 0 \}$. In the case $d>0$, $D_C$ is called non-degenerate if 
\begin{equation} \label{delbar-generic}
D_C\oplus ev_\Omega : C^\infty (N_C) \rightarrow C^\infty ( N_C \otimes T^{0,1}C) \oplus  _{p\in \Omega} N_p 
\end{equation}
has trivial cokernel. 

\begin{definition} \label{def-generic}
A pair 
$(J,\Omega)\in \mathcal{A} _m (X) $, $m\ge 0$,  is said to be generic if the following
five conditions are met for all $E\in H_2(X;\mathbb{Z})$ for which the number $d$ as defined by \eqref{indexofd} is no greater than $m$
(see \cite{kn:taub4} for more details, especially on the  definition of $n$-non-degeneracy which is immaterial for the present discussion and 
thus omitted):

\begin{enumerate}
\item For a fixed class $E\in H_2(X;\mathbb{Z})$, there are only finitely many embedded $J$-holomorphic curves representing $E$ and 
containing $d$ points of $\Omega$. 
\item For each $J$-holomorphic curve $C$, the operator $D_C$ is non-degenerate. 
\item There are no connected $J$-holomorphic curves representing the class $E\in H_2(X;\mathbb{Z})$ and containing more than $d$ points of $\Omega$. 
\item There is an open neighborhood of $(J,\Omega)$ in $\mathcal{A} _m (X)$ such that each pair $(J',\Omega ')$ from that neighborhood satisfies conditions 1-3 above. 
Furthermore, the number of $J'$-holomorphic curves containing $d$ points of $\Omega '$ is constant as $(J',\Omega ')$ varies through the said neighborhood. 
\item  If $E^2=K\cdot E = 0 $ then each of the finitely many $J$-holomorphic curves in $E$ containing $d$ points of $\Omega$ is $n$-non-degenerate for each 
positive integer $n$.  
\end{enumerate}
\end{definition}

The set of generic pairs $(J,\Omega)$, which we denote by  $\mathcal{J} ^{reg}_d (X) $ (or simply by $\mathcal{J} ^{reg} (X) $  when no confusion is possible), is a Baire subset of $\mathcal{A} _d (X) $. 

We are now ready to state our main result:
\begin{theorem} \label{theo-main}
Let $(X,\omega)$ be a minimal symplectic 4-manifold and $C$ a connected, embedded symplectic surface in $X$ of genus $g\ge 1$ and 
with $C^2 \ge g-1$. 
Then for any $\delta>0$ there exists a generic pair $(J_\delta , \Omega_\delta)  \in \mathcal{J} ^{reg} (X) $  and a connected $J_\delta$-holomorphic curve $C_\delta$ 
inside the radius $\delta$
tubular neighborhood of $C$, isotopic to $C$. Furthermore, $C_\delta$ contains all $d$ points of $\Omega_\delta$. 
\label{main}
\end{theorem} 
\begin{corollary} \label{corollary1}
The above theorem remains true if $C=\sqcup \, C_i$ is a disjoint union of connected symplectic manifolds provided the condition $C_i^2 \ge g_i - 1$ holds for 
each component $C_i$. That is, one can find a curve $C_{\delta} = \sqcup \,  C_{\delta , i } $ where each $C_{\delta, i}$ is an isotopic translate of $C_i$ 
inside a radius $\delta$  tubular neighborhood of $C_i$. 
\end{corollary}
\begin{remark} \label{negsquare1}
Whenever $(J,\Omega)$ is a generic pair in the sense of definition \ref{def-generic}, the Gromov-Witten moduli space $\mathcal{M}^{Gr}_X(E)$ is a smooth 
manifold of (real) dimension $2d \ge 0$. This together with the adjunction formula for a connected $J$-holomorphic curve 
$C\in \mathcal{M}^{Gr}_X(E)$ implies that $E^2 \ge g-1$ (where $g$ is the genus of $C$). Conversely, given a connected symplectic curve $C$ of 
genus $g$ satisfying $C^2\ge g-1$, theorem \ref{theo-main} shows that there are no other obstructions for the existence of a generic pair $(J,\Omega)$ 
making $C$ into a $J$-holomorphic curve.  
\end{remark}
\begin{remark} \label{negsquare2}
Suppose that $(J, \Omega)$ is a generic pair and let $C$ be a connected $J$-holomorphic curve of genus $g$ and with $[C]=E$. The inequality $E^2 \ge g - 1$ 
from the previous remark, shows that $J$-holomorphic curves with negative square can only occur when  $E^2= -1$ and $g=0$. This case however is 
excluded if $X$ is a minimal manifold (as is assumed in theorem \ref{theo-main}).   
\end{remark}

It is interesting to compare the result of theorem \ref{main} to the result proved in \cite{hofer}. Expressed in our notation, among 
other results, it is proved in \cite{hofer} that for $C^2 \ge 2g -1 $, the operator $D_C$ is surjective for {\bf any} choice of an 
almost-complex structure $J$ compatible with the symplectic form $\omega$. The improvement of the inequality in 
theorem \ref{main} comes at the twofold expense of first not 
being able to choose the almost-complex structure arbitrarily but rather 
from a dense (second-category) subset of almost-complex structures. Secondly, one may have to slightly 
\lq\lq wiggle\rq\rq \, $C$ to get the desired curve.  We would also like to remark that the case 
of genus 0, which is excluded from theorem \ref{main}, is completely covered by the results of \cite{hofer}. 

The proof of theorem \ref{main} rests on the observation that the property of a $J$-holomorphic curve $C$ to be generic with respect to a pair 
$(J,\Omega)\in \mathcal{A}_d (X)$ is local in nature, that is, it only depends on the 
restriction of the of $(J,\Omega)$ to a tubular neighborhood $N(C)$ of the curve $C$.  By the 
symplectic neighborhood theorem for four-manifolds (cf. \cite{dusa}), $N(C)$ is up to symplectomorphism determined by its volume and the square $C^2$ of the curve
$C$.  
Thus one is led to search for universal models of symplectic four-manifolds $Y_{g,n}$ with a Gromov-Witten basic class $E_{g,n} \in H_2(Y_{g,n};\mathbb{Z})$ with 
$E_{g,n}\cdot E_{g,n} = n$ and for which a connected genus $g$ $J$-holomorphic representative exists for all generic $(J,\Omega)$. These manifolds together with 
their Gromov-Witten invariants are discussed in section \ref{three-two} after a brief survey of Seiberg-Witten theory on four-manifolds with $b^+ = 1$ which is given 
in section \ref{three-one}. No originality is claimed on any of the facts stated in section \ref{three}, they serve merely as a reminder and 
to set notation. The proof of theorem \ref{main} is then completed in section \ref{four}. Section \ref{two} gives applications of the main theorem. 
\section{Applications} \label{two}

As an application of theorem \ref{main}, we give examples of symplectic manifolds with non-empty Seiberg-Witten and Gromov-Witten moduli spaces under 
generic conditions, whose associated 
invariants are zero. Such examples can be found for the case where the dimension of the moduli space is  zero as well as for the case of positive dimension. 
\vskip2mm

{\bf Example 1: } Consider the elliptic surface $E(n)$. It has a symplectic section $S_n$ with genus zero and square $-n$. Let $F_i$, $i=1,2,...$ be regular fibers of 
the elliptic fibration. Then the symplectic surface $C_{n,m}$, obtained by smoothing the surface $S_n \cup F_1 \cup ... 
\cup F_m$, is a genus 
$g_{n,m}=m$ surface of square $2m-n$. Choosing
$m\ge n-1$ ensures the condition $C_{n,m}^2\ge g_{n,m} -1$. Theorem \ref{main} provides a generic pair $(J,\Omega) \in \mathcal{J} ^{reg} (E(n)) $ and a $J$-holomorphic 
curve $C'_{n,m}$ in the class $[C_{n,m}]$. 
In particular, the moduli space $\mathcal{M} ^{Gr}_{E(n)} ([C_{n,m}])$ for this generic pair $(J,\Omega)$ is nonempty while $Gr_{E(n)}([C_{n,m}]) = 0 $ . 
The dimension of the moduli space is $$dim _\mathbb{R}\, \mathcal{M} ^{Gr}_{E(n)} ([C_{n,m}]) = 2(m-n+1)$$ 
\vskip2mm

{\bf Example 2: } Let $\Sigma$ be a genus 2 Riemann surface and let $X=\Sigma\times T^2$. Choose the symplectic form 
$\omega$ on $X$ to be the sum of volume forms $\omega _\Sigma$ and $\omega _{T^2}$ 
on $\Sigma$ and $T^2$  for which Vol$(\Sigma)$=1=Vol$(T^2)$. Let $C$ be the symplectic surface obtained by smoothing 
$\Sigma \cup  T^2 $. Then 
the genus of $C$ is 3 and its square is 2, in particular, dim$\mathcal{M} _X^{Gr}([C])=0$ and dim$\mathcal{M} _X^{SW}(L)=0$ for $L=2P.D.([C])-K$. 

Pick an almost-complex structure  $J\in \mathcal{J} ^{reg}(X) $ ($\Omega$ is just the empty set here and we suppress it from the notation) 
and a $J$-holomorphic curve $C'$ in the class $[C]$. 
It is not hard to see, but somewhat tedious, that all $J$-holomorphic curves in $[C]$ are connected curves of genus 3. To see this, consider the two possible 
alternatives: 
\begin{enumerate}
\item There is a representative $D'$ of $[C]$ of the form $D' = D'_1\sqcup ... \sqcup D'_n$ with $D'_i\cdot D'_i=1$ for $i=1,2$ and $D'_i\cdot D'_i=0$ for $i\ge 3$. 
This is an immediate contradiction since classes of square 1 cannot exist on a manifold with even canonical class.
\item There is a representative $D$ of $[C]$ of the form $D = D_1\sqcup ... \sqcup D_n$ with 
$D_1^2  = 2$ and $D_i^2 = 0 $ for $i\ge 2$. 
This implies that $g(D_i) = 0 $ for $i\ge 2$ and $2\le g(D_1) \le 3$. The latter claim follows readily from the fact that the dimension $dim\, \mathcal{M}_X^{Gr} ([D_1]) = 2(D_1^2 - g(D_1) +1) $ is non-negative and from the adjunction formula for $D$. 
The case $g(D_1)=2$ leads (via the adjunction formula applied to $[C]$) to $[C]\cdot K = 0$ 
which is a contradiction. Thus the only possibility is $g(D_1) = 3$ implying $K\cdot D_1 = 2$. 

Since $\omega \in H^2(X;\mathbb{Z})$ and $\omega ([C])=2$,  
we see immediately that $n\le 2$. Suppose thus that $D=D_1 \sqcup D_2$. Then by $K\cdot D_1=2 $ we see that 
$[D_1] = [\Sigma] + a [T^2] + F$ where $F\in H_2(X;\mathbb{Z})$ is generated by classes obtained from cross-products of 
1-cylces on $\Sigma$ with 1-cylces on $T^2$. This forces 
$[D_2] = (1-a) [T^2] -F$.
Notice that $F\cdot \Sigma = F\cdot T^2 = \omega \cdot F = 0$.  
From $D_1^2 = 2$ we infer that $2a+F^2 = 2$ and from $D_2^2 = 0$ we get $F^2=0$. Thus $a=1$ and so 
$[D_1] = [\Sigma] + [T^2] + F$ and $[D_2] = - F$. This now leads to a contradiction since now $\omega (D_2) =0$ and so
$D_2$ cannot be a $J$-holomorphic curve. 
\end{enumerate}
Each point in $\mathcal{M} ^{Gr}_X([C])$ gives rise to a Seiberg-Witten monopole in $\mathcal{M}^{SW}_X(L)$ with $L=2\,P.D.([\Sigma])$
(see \cite{kn:taub3}). It was shown in \cite{kn:me}
that each such monopole is a smooth point in the moduli space for large enough values of $r$ in the Taubes perturbation form $\mu_0 = F_0^+ -ir\omega/8$. Said in other words, 
the pair of metric and perturbation form $(g,\mu_0)$ (with $g$ being the metric induced by $\omega$ and $J$) is a generic pair for the Seiberg-Witten theory for the 
\spin $L$ and as such gives rise to a smooth moduli space.  On the other hand $SW_X(L)=0$ as can be seen in a number of 
ways (for example, introduce the \lq\lq twisted\rq\rq  
symplectic form $\omega ' = 1.1 \omega _\Sigma + \omega _{T^2}$. Then $L\cdot \omega ' > K\cdot \omega '$ 
which according to \cite{taubes} implies that $L$ cannot be a basic class). 
\section{Preliminaries} \label{three}
\subsection{Seiberg-Witten theory on manifolds with $b^+=1$} \label{three-one}

Let $X$ be a 4-manifold with $b^+=1$. For a given $Spin^c$-structure $W=W^+\oplus W^-$, with determinant $L=$det$(W^+)\in H^2(X;\mathbb{Z})$, the Seiberg-Witten
invariant depends on a choice of a chamber inside the space $Met \times i \, \Omega ^{2,+}$. Here $Met$ is the space of Riemannian metrics on $X$. 
The two chambers are divided by a (real) codimension 1 wall of pairs $(g,\mu )$, defined by the equation 
$$  \frac{i \mu }{2\pi} \wedge \omega _g  - L\wedge \omega _g = 0$$
where $\omega _g$ is a generator of the positive forward cone in $H^2(X;\mathbb{Z})$. In the case where $X$ is symplectic, we agree to always choose $\omega _g$ 
to be the symplectic form.

The Seiberg-Witten equations do not  admit reducible solutions if $(g,\mu )$ doesn't  lie on the wall. We denote the two chambers by $\mathcal{C}^- (L)$ and 
$ \mathcal{C}^+  (L)$ according to the sign of the expression
$$ \langle  \frac{i \mu }{2\pi} \wedge \omega _g    -    L\wedge \omega _g , [X]                        \rangle   $$
We will denote the Seiberg-Witten invariant by $SW_X^{\pm} (L)$ according to the choice of chamber $\mathcal{C} ^{\pm}(L)$ from which the pair $(g,\mu )$ 
used in calculating the invariant, was taken from. The number $SW_X^+(L) - SW_X^- (L)$ is called the {\sl wall crossing number} and it is well understood (see for example 
\cite{kn:li}). The 
special case relevant to the present situation is stated in the following theorem (Corollary 1.4 in  \cite{kn:li}):
\begin{theorem} \label{theorem:wall} Let $X$ be an $S^2$-bundle over a Riemann surface $\Sigma$ of genus $g$.
Let $E\in H^2(X;\mathbb{Z})$ with 
$(2E+c_1(X))^2 \ge 2e_X + 3\sigma _X$. Then the wall crossing number is 
$$ SW_X^+(L) - SW_X^- (L) = \pm \left( \frac{2E+ c_1(X)}{2}\, [S^2]   \right)^g $$ 
where $[S^2]$ is the fiber class. 
\end{theorem} 
%
%
%
\subsection{The Gromov-Witten Invariants of $Y_0$ and $Y_1$} \label{three-two}

This section describes the spaces $Y_{g,n}$ mentioned in the introduction as well as their Gromov-Witten basic classes $E_{g,n}$. 
As it turns out, it suffices to consider only two symplectic manifolds $Y_0$ and $Y_1$ by letting $Y_{g,2n} = Y_0$ and $Y_{g,2n-1} = Y_1$. 
The main results of this section, corollaries 
\ref{firstcoro} and \ref{secondcoro}, are well known and their proofs can be found in the literature (see e.g. \cite{friedman}). 
They are only 
included here for continuity of argument and for the benefit of the reader,  no originality is claimed. 

Let $\Sigma$ be any Riemann surface of genus $g\ge 2$. Define the $Y_0$ and $Y_1$ to be
$$ Y_0 = \Sigma \times S^2 \quad \quad \mbox{and} \quad \quad Y_1= Y_0 \# _{F_0 = F_1} ( S^2 \tilde{\times} S^2 ) $$ 
In the above, $S^2 \tilde{\times} S^2$ denotes the twisted $S^2$ bundle over $S^2$. It is diffeomorphic to $\mathbb{CP}^2\# \overline{\mathbb{CP}^2}$. 
As $Y_1$ is obtained by fiber sum of  two $S^2$ fibrations, it itself inherits the 
structure of an $S^2$ fibration over $S^2$. 

To calculate the Gromov-Witten invariants of $Y_i$, we invoke Taubes' theorem relating the Gromov-Witten invariants to the Seiberg-Witten invariants, 
the latter of which often prove easier to calculate. The following theorem can be found in \cite{kn:taub3}.
\begin{theorem} \label{taubes}
Let $(X,\omega )$ be a symplectic 4-manifold with $b^+=1$. Let $\mu_0 = F^+_{A_0} - ir\omega /8\in i\Omega ^{2,+}$ (where $A_0$ is a certain connection on 
the canonical line bundle) and let 
$g$ be any generic metric compatible with the symplectic form. Then, for any $E\in H^2(X;\mathbb{Z})$, the 
Seiberg-Witten invariant of $X$ for the \spin $W^+ _E=E\oplus (E\otimes K^{-1})$, calculated with the metric $g$ and perturbation form $\mu_0$ with $r\gg 1$, is equal 
to the Gromov-Witten invariant for the class $E$. 
\end{theorem}  
The Seiberg-Witten invariants for both $Y_0$ and $Y_1$ are calculated in much the same way. We will explicitly only give  
the calculation for $Y_0$ here and refer to the minute differences that occur for $Y_1$.

The main input for calculating the Seiberg-Witten invariants of $Y_0$ and $Y_1$ are the wall crossing formula and the existence of metrics with positive scalar curvature. 

Let $g_\Sigma$ and $g_{S^2}$ be metrics on $\Sigma$ and $S^2$ with constant scalar curvature and with volumes equal to $4\pi (g-1)$ and $4\pi$ 
respectively. It follows from the Gauss-Bonnet theorem that the scalar curvatures $s_\Sigma$ and $s_{S^2}$ of these metrics are 
$$ s_\Sigma = -1 \quad \quad \mbox{ and } \quad \quad  s_{S^2} = 1 $$
Denote by $\omega _\Sigma$ and $\omega _{S^2}$ the volume forms induced by $g_\Sigma$ and $g_{S^2}$ and define the symplectic form  
$\omega_{\lambda ,\varepsilon}$ on $Y_0$ to be 
\begin{equation}
\omega_{\lambda , \varepsilon}  = \lambda \cdot \omega_{\Sigma} + \varepsilon\cdot\omega _{S^2} \label{eq:sympform}
\end{equation}
The positive 
parameters $\lambda, \varepsilon >0$ will be chosen later, $\varepsilon$  should be thought of as being small.  
The 
product metric 
$$g_{\lambda , \varepsilon}  = {\lambda} \, g_\Sigma \oplus {\varepsilon} \, g_{S^2} $$  
on $Y_0$ is compatible with $\omega _{\lambda , \varepsilon}$
and 
its  scalar curvature $s_{\lambda , \varepsilon}$ is
$$ s_{\lambda , \varepsilon} = - \frac{1}{\lambda} + \frac{1}{\varepsilon }  $$  
Our first condition on the parameters $\lambda$ and $\varepsilon$ will be that $\varepsilon < \lambda$, ensuring that $s_{\lambda , \varepsilon} > 0$ (the choice of the 
second condition is deferred to section \ref{four}). 

With $\omega _{\lambda , \varepsilon} $ chosen as in \eqref{eq:sympform}, the canonical class $K_0$ of $Y_0$ is easily calculated from the adjunction formula and 
from the fact that both $\Sigma \times \{ pt\}$ and $\{pt \} \times S^2$ are symplectic submanifolds of $Y_0$. One finds that
$$ K_0 = (2g-2)\, \overline{S} - 2 \, \overline{\Sigma}\,  \in \, H^2(Y_0;\mathbb{Z}) $$
where $\overline{S} = P.D. ([S^2])$ and $\overline{\Sigma} = P.D. ([\Sigma])$. 

We will label $Spin^c$-structures of $Y_0$ by elements $E\in H^2(Y_0;\mathbb{Z})$ by letting $W_E$ be the $Spin^c$-structure with 
$W^+ _E = E \oplus (E\otimes K^{-1})$.
Thus the determinant line bundle $L=\, $det$(W_E^+)$ is equal to $2E-K$. We label the corresponding 
Seiberg-Witten moduli spaces by $\mathcal{M} _{Y_0}^{\pm} (L)$, the sign again depending upon the chamber $\mathcal{C} ^{\pm}(L)$ determined by the metric and perturbation.

For $a,b\in \mathbb{Z}$, let $E=a\, \overline{\Sigma} + b\, \overline{S}$ and consider the $Spin^c$-structure $W_E$. The dimension of the Seiberg-Witten 
moduli space is given by 
$$ \mbox{dim}_{\mathbb{R}}\, \mathcal{M} ^{\pm}(E) = \frac{1}{4} \left( L^2  - K_0^2 \right) = 2b\, (a+1) - a\, (2g-2)$$
In order for the $Spin^c$-structure $W_E$ to have nonzero Seiberg-Witten invariant, the dimension of the moduli space needs to be non-negative. 
In the case of $a=1$ (the case of interest to us) together with the  observation that $E^2 = 2b$, the above formula leads to a 
necessary condition for the nonvanishing of the 
invariant:
$$ E^2 \ge g-1$$
Consider now $E=\overline{\Sigma} + b\, \overline{S}$ with $E^2 = 2b\ge g-1$ and let $L=2E-K$. It is easy to see that  
\begin{equation}
\langle L\wedge \omega , [Y_0] \rangle = 32\pi \lambda (g-1) + 16\pi \varepsilon (b-g+1)   \label{eq:secondcond}
\end{equation}
Two pairs of a metrics and perturbation forms will play a role in the subsequent discussion: 
\begin{enumerate} 
\item $(g,\mu) = (g_{\lambda, \varepsilon},0)$:  By our choice $\lambda > \varepsilon$ and by the restriction $2b\ge g-1$, the right-hand side of \eqref{eq:secondcond} is positive:
$$32\pi \lambda (g-1) + 16\pi \varepsilon (b-g+1) \ge 16\pi (g-1)(2\lambda - \frac{1}{2} \varepsilon ) > 0 $$
This means that the pair $(g_{\lambda, \varepsilon}, 0 )$ lies in the chamber $\mathcal{C} ^-(L)$.
\item $(g,\mu) = (g_0,\mu_0)$:  Here $g_0$ is any generic metric (but still compatible with $\omega _{\lambda, \varepsilon}$) and $\mu_0$ is 
Taubes' perturbation form 
$$\mu_0 = F_{A_0} ^+ - \frac{ir\omega }{8}$$
It is easily checked that for large enough $r$, the pair $(g_0,\mu _0)$ lies in $\mathcal{C} ^+(L)$ (for any $Spin^c$-structure). 
\end{enumerate} 

By the positivity of $s_{\lambda, \varepsilon}$ we have that $SW^-_{Y_0} (L) = 0 $ which together with 
theorems \ref{theorem:wall} and \ref{taubes} immediately gives 
\begin{corollary} \label{firstcoro}
For $g\ge 1$, let $E_{g,2n}  = \overline{\Sigma} + n\, \overline{S^2} \in H^2(Y_0;\mathbb{Z})$ with $E_{g,2n} ^2\ge g-1$. Then 
$$ Gr_{Y_0} (E_{g,2n}) = \pm \, 2^g  $$
\end{corollary} 
While the discussion preceding corollary \ref{firstcoro} was  for  the case $g\ge 2$, it is not hard to see that it still remains valid in 
the case $g=1$. The changes that need to be made to the analysis preceding the corollary are: 
choose the product metric on $\Sigma = T^2$ so that 
its scalar curvature is zero. Choose $\omega_{\lambda , \varepsilon}$ and $g_{\lambda , \varepsilon}$ as before and 
observe that $s_{\lambda , \varepsilon} = 1/\varepsilon$ which is positive for $\varepsilon > 0$. 
The rest of the discussion goes over verbatim and 
so establishes the validity of corollary \ref{firstcoro} in the case $g=1$ as well. 

We finish this section by showing that an analogous result holds for $Y_1$. 
In $Y_1$, let $\Sigma' = \Sigma {\#} S\subseteq Y_0 \# _{F_0=F_1}   (S^2 \tilde{\times} S^2)$  with 
$S=\mathbb{CP}^1\subseteq \mathbb{CP}^2\# \overline{\mathbb{CP}^2} \cong S^2 \tilde{\times} S^2$. Let $F$ denote 
a fiber of the fibration $Y_1\rightarrow S^2$. The canonical class $K_1$ of $Y_1$ is 
$$ K_1 =  (2g-1)\,\overline{F} - 2\, \overline{\Sigma}' _0 \quad \quad \overline{F} = P. D.([F]),\, \overline{\Sigma}' = P.D.([\Sigma'])   $$ 
As with $Y_0$, consider $E=a\overline{\Sigma}' _0 + b \overline{F} \in H^2(Y_1;\mathbb{Z})$. The dimension for the 
Seiberg-Witten moduli space for the \spin $W_E$ is 
$$ \mbox{dim}\mathcal{M} ^{SW}_{Y_1}(L) = 2b(a+1)  - 2a (g-2)  $$ 
In the case when $a=1$,  the  necessary condition for the nonvanishing of $SW^\pm _{Y_1} (L)$  (with $L=2E-K_1$) becomes 
$$ E^2 = 2b+1 \ge g-1 $$
It is a known fact (cf. \cite{kn:lebrun}) that ruled surfaces admit metrics of positive scalar curvature. The rest of the discussion for $Y_1$ proceeds now in much the 
same way as that for $Y_0$ and one arrives at the following analogue of corollary \ref{firstcoro}: 
\begin{corollary} \label{secondcoro}
Let $E_{g,2n+1} = E= \overline{\Sigma_0} + n\, \overline{F} \in  H^2(Y_1;\mathbb{Z})$ with $E^2\ge g-1$. Then 
$$ Gr_{Y_1} (E) = \pm \, 2^g  $$
\end{corollary} 
%
%
%
\section{ Proof of Theorem \ref{main}} \label{four} 

We now proceed to the proof of the theorem \ref{main}. Let $C$ be an embedded, connected, symplectic submanifold of $(X^4, \omega)$ of genus $g\ge 1$ and with 
square $[C] ^2 = n \ge g-1$. Assume in addition that $n=2k$ is even, the case where $n$ is odd is treated in much the same way by 
replacing $Y_0$ below with $Y_1$. 
Let $N(C)$ be a tubular neighborhood of $C$ in $X$ and let Vol$(C)$ be the volume 
of $C$. 

On the other hand, let $D$ be any of the (at least) $2^g$ $J'$-holomorphic curves in $Y_0$ in the class $[\Sigma] + k [S^2]$ for the choice of a generic  
pair $(J',\Omega ')$ on $Y_0$. This last statement uses corollary \ref{firstcoro} (or corollary \ref{secondcoro} in the case of $n=2k-1$). 
Adjust the choices of $\lambda$ and $\varepsilon$ so that the Vol$(D) = $Vol$(C)$ (in addition to $\lambda > \varepsilon >0$). 
Let $N(D)$ be a tubular neighborhood of $D$ in $Y_0$ containing no other $J'$-holomorphic curves besides $D$. 

By the symplectic neighborhood theorem for 4-manifolds (cf. \cite{dusa}, exercise 3.30), the tubular neighborhood of a connected,  embedded symplectic surface is up to 
symplectomorphism determined by the square and volume of the surface. We would like to say that the pairs 
$(N(C),\omega|_{N(C)})$ and $ (N(D), \omega_{\lambda, \varepsilon}|_{N(D)})$ are symplectomorphic via a symplectomorphism 
$\varphi : N(C) \rightarrow N(D)$ taking $C$ to $D$. There is one potential problem with this approach and that is that a priori all of the at least $2^g$ $J'$-holomorphic 
curves in the class $[\Sigma] + k [S^2]$ in $Y_0$ may be disconnected. Fortunately, the opposite extreme is true as the next lemma contests:
\begin{lemma}
Let $(J',\Omega ')$ be a generic pair on $Y_i$ and let $D$ be an embedded $J'$-holomorphic curve in $Y_i$ containing $\Omega '$.  
Suppose that $D$ represents the homology class $[\Sigma] + k [S^2]$ in the case $i=0$ and 
represents the class $[\Sigma '] + k [F]$ in the case $i=1$. Then $D$ is connected. 
\end{lemma}
\begin{proof}  Assume to the contrary that we can write $D$ as a disjoint union 
$D=D_1 \sqcup D_2$. We will show that one of the two components has fundamental class zero. 

{\bf Case of $i=0$: }
Let $[D_1] = a [\Sigma] + b[S^2]$ and $D_2 = c[\Sigma] + d[S^2]$. Since $a+c=1$ we 
can assume that $a\ge 1$. We will first show that in fact $a=1$ and thus $c=0$. 

It is a well known fact that for generic almost-complex structures, $J$-holomorphic curves intersect non-negatively (see \cite{dusazero}). Observe also that 
the manifolds $Y_i$ are minimal and so remark \ref{negsquare2} applies (excluding the existence of $J'$-holomorphic curves with negative square). 
We know by 
corollary \ref{firstcoro} that for $N$ large enough, the class $[\Sigma] + N [S^2]$ has $J$-holomorphic representatives. Thus we get 
\begin{align}
 [D_2] \cdot ( [\Sigma] + N [S^2] ) \ge 0 \quad & \Longrightarrow  \quad c\, N + d \ge 0 \cr
						& \Longrightarrow \quad (1-a) N + d \ge 0 \cr
						& \Longrightarrow \quad 1 + \frac{d}{N} \ge a \ge 1 \cr
						& \Longrightarrow \quad a=1 \mbox{ and } c=0 
\end{align}
Since $D_1$ and $D_2$ are disjoint, we find that $0 = [D_1]\cdot [D_2] = d$ which shows that $[D_2]= 0$. 

{\bf Case of $i=1$: } Let $[D_1] = a [\Sigma ' ] + b[F]$ and $D_2 = c[\Sigma '] + d[F]$. Since as before we have $a+c=1$ we 
can again assume that $a\ge 1$. Using corollary \ref{secondcoro} we know that the class $[\Sigma' ] + N [F]$ has $J$-holomorphic representatives for all 
sufficiently large N. Then arguing as above we have:
\begin{align}
 [D_2] \cdot ( [\Sigma ' ] + N [F] ) \ge 0 \quad & \Longrightarrow  \quad c+ c\, N + d \ge 0 \cr
						& \Longrightarrow \quad (1-a) (N+1)  + d \ge 0 \cr
						& \Longrightarrow \quad 1 + \frac{d}{N+1 } \ge a \ge 1 \cr
						& \Longrightarrow \quad a=1 \mbox{ and } c=0 
\end{align}
The fact $0=[D_1]\cdot [D_2] = d$ completes the proof. 
\end{proof}

Use $\varphi$ together with $J'$ on $N(D)$ to induce an almost-complex structure (still denoted by $J'$) on $N(C)$.  Extend $J'$ over all of $X$ 
in an arbitrary manner and denote it by $J''$. Let $\Omega ''$ denote the set $\varphi ^{-1} (\Omega ')$.   

Observe that $(J'',\Omega '') \in \mathcal{J} ^{reg}_d (N(C))  $ but it could happen that $(J'',\Omega '') \notin \mathcal{J} ^{reg}_d (X)  $ as there may be other $J''$-holomorphic 
curves in $X$ for which the operator defined in \eqref{delbar-generic} is not surjective. 
However, generic pairs  $(J,\Omega) $ on $X$ 
are dense in $\mathcal{A} _d(X) $ and so we can find, in an arbitrarily small neighborhood of $(J'',\Omega '')$, a pair $(J,\Omega)$ that is generic.  The following standard 
proposition completes the proof of theorem \ref{main}. 
\begin{prop} \label{lastprop}
Let $\varepsilon >0$ be arbitrary. Then there exists $\delta > 0$ such that if $$\mbox{dist }[ (J,\Omega), (J'',\Omega '') ]< \delta $$ then there exists a $J$-holomorphic curve $C'$ 
in an $\varepsilon$ tubular neighborhood of $C$. 
\end{prop}  
\begin{proof} This is a direct consequence of the fourth point in the definition of genericity applied to the two pairs $(J'' , \Omega '' )|_{N(C)}$ and 
$(J , \Omega  )|_{N(C)}$. By construction $(J'',\Omega '') \in \mathcal{J} ^{reg} _d (N(C))$ and  clearly 
$$\mbox{dist }[ (J,\Omega)|_{N(C)} , (J'',\Omega '')|_{N(C)}  ] \le \mbox{dist }[ (J,\Omega), (J'',\Omega '') ]    $$
This completes the proof of the proposition as well as theorem \ref{main}. 
\end{proof}
\vskip3mm

\begin{proof}[Proof of corollary \ref{corollary1}]  The proof of corollary \ref{corollary1} proceeds in much the same way. For each component $C_i$ of $C$, one 
finds a generic pair $(J'_i, \Omega _i ')$ on a tubular neighborhood $N(C_i)$ of $C_i$. One extends the almost-complex structers $J_i'$ to an arbitrary 
almost-complex structure $J''$ on $X$ and defines $\Omega '' = \sqcup \Omega _i ''$ where the $\Omega _i ''$ are defined as $\Omega ''$ was  in the proof of theorem \ref{main}. The analogue of propositon \ref{lastprop} completes the proof of corollary \ref{corollary1}. 
\end{proof}

\end{document}